\documentclass[12pt]{article}

\usepackage[latin1]{inputenc}\usepackage[english]{babel}
\usepackage{amsmath}
\usepackage{amsfonts}
\usepackage{latexsym}
\usepackage[dvips]{graphicx}
\usepackage{url}

\newtheorem{theorem}{Theorem}
\newtheorem{rem}[theorem]{Remark}

\newtheorem{lemma}[theorem]{Lemma}

\newtheorem{eg}[theorem]{Example}

\newtheorem{conj}[theorem]{Conjecture}

\newtheorem{prob}[theorem]{Problem}

\newcommand{\pf}{\noindent {\bf Proof: }}

\begin{document}
\title{Recurrent words with constant Abelian complexity}
\author{James Currie\thanks{The author is
supported by an NSERC Discovery Grant.} and
Narad Rampersad\thanks{The author is supported by an NSERC Postdoctoral
Fellowship.} \\
Department of Mathematics and Statistics \\
University of Winnipeg \\
515 Portage Avenue \\
Winnipeg, Manitoba R3B 2E9 (Canada) \\
\url{j.currie@uwinnipeg.ca} \\
\url{n.rampersad@uwinnipeg.ca}}

\maketitle

\begin{abstract}
\noindent We prove the non-existence of recurrent words with constant
Abelian complexity containing 4 or more distinct letters. This answers
a question of Richomme et al.\vspace{.1in}
\end{abstract}

\section{Introduction}
One of the central notions in combinatorics on words is that of the
subword complexity of an infinite word.  Richomme, Saari, and Zamboni
\cite{RSZ09} have recently begun a systematic study of the Abelian
analogue of the subword complexity of infinite words.  In this paper
we resolve one of the open problems from their study by showing the
non-existence of recurrent words with constant Abelian complexity
containing 4 or more distinct letters.

Let $\Sigma$ be a finite alphabet and let $\Sigma^*$ be the set of all
finite words over the alphabet $\Sigma$.
Consider the equivalence relation $\sim$ on  $\Sigma^*$, defined
by $$u\sim v\mbox{ if }u\mbox{ is an anagram of }v.$$ Thus $1232\sim
2132$. We write $[u]$ for the equivalence class of
$u$ under $\sim$. For example, $[121]=\{112, 121, 211\}$. We call
$[u]$ an {\bf Abelian word}. If $u$ is a factor of a  word $w$, we call $[u]$
an {\bf Abelian factor} of $w$.  The length of an Abelian factor is
the length of any one of its representatives.

If $w$ is an infinite word, the {\bf subword complexity function} of
$w$ is the function $f : \mathbb{N} \to \mathbb{N}$, where for $m =
1,2,\ldots$, the value of $f(m)$ is the number of factors of $w$ of
length $m$.  Similarly, the {\bf Abelian complexity function} of $w$ is the
function $\tilde{f}(m) : \mathbb{N} \to \mathbb{N}$, where for $m =
1,2,\ldots$, the value of $\tilde{f}(m)$ is the number of Abelian
factors of $w$ of length $m$.

An infinite word $w = w_0w_1\cdots$, where $w_i \in \Sigma$ for $i =
0,1,2\ldots$, is {\bf ultimately periodic} if there exist a
non-negative integer $c$ and a positive integer $p$ such
that $w_i = w_{i+p}$ for all $i \geq c$.  A classical result of
Morse and Hedlund \cite{MH40} shows that an infinite word $w$ is ultimately
periodic if and only if its complexity function $f$ is eventually
constant.  If $w$ is not ultimately periodic, then $f(m) \geq m+1$ for
all $m$.  The well-studied {\bf Sturmian words} are precisely the
aperiodic words of minimal complexity (i.e., those words for which
$f(m) = m+1$ for all $m \geq 1$).  Coven and Hedlund \cite{CH73}
showed that any Sturmian word has constant Abelian complexity.  In
particular, for any Sturmian word, one has $\tilde{f}(m) = 2$ for all
$m \geq 1$.

Sturmian words are necessarily over a binary alphabet; it is therefore
natural to ask if over an $n$-letter alphabet, where $n \geq 3$,
there is an infinite word $w$ with Abelian complexity function
$\tilde{f}(m) = n$ for all $m \geq 1$.  Without further qualification,
this question is not very interesting, as one easily sees that the word
$$123\cdots (n-1)nnnnnnnn\cdots$$
over the alphabet $\{1,2,\ldots,n\}$ has exactly $n$ Abelian factors
of each length $m \geq 1$.

This observation leads us to the following definition.  We say that a
word $w$ is {\bf recurrent} if every factor of $w$ occurs infinitely
often in $w$.  Any Sturmian word is recurrent, so such words provide
examples of recurrent words with constant Abelian complexity over a
binary alphabet.  Richomme, Saari, and Zamboni showed that there are
recurrent words over a $3$-letter alphabet with exactly $3$ Abelian
factors of each length $m \geq 1$, thereby answering a question of
Rauzy.  They also posed a question of their own, namely, ``Does there
exist a recurrent word over a $4$-letter alphabet with exactly $4$
Abelian factors of each length?''  They also conjectured that the
answer to the question should be ``no''.  We show that this is indeed
the case.  Moreover, our main result also applies to alphabets of size
greater than $4$.

\begin{theorem}\label{main}
Let $n \geq 4$ be an integer.  There is no recurrent word over an
$n$-letter alphabet with exactly $n$ Abelian factors of each length
$\geq 1$.
\end{theorem}

\section{Proof of Theorem~\ref{main}}
Fix a positive integer $n\ge 4$. Let $\Sigma$ be the alphabet
$\{1, 2, 3, \ldots, n\}$. Let $w$ be a finite or infinite word.
Consider the graph $G$ with vertex set $\Sigma$, and an edge $ij$
whenever at least one of $ij$ and $ji$ is a factor of $w$. Note
that $G$ may contain loops, but not multiple edges.\footnote{For
definiteness of notation, let us say that we never call $a$ a
neighbour of itself; however, we will count a loop based at $a$ as
contributing 1 to the degree of $a$. Thus the degree of a vertex
$a$ in $G$ will be the number of distinct neighbours of $a$, plus
the number of loops based at $a$; since we do not allow multiple
edges, the number of loops based at $a$ is 0 or 1.} From now on
suppose that $w$ is a fixed recurrent word, having constant
Abelian complexity $n$.

\begin{lemma}\label{one cycle}Graph $G$ consists of a spanning tree and one additional edge. Thus $G$ contains a unique cycle $C$ (which is possibly a loop).
\end{lemma}
\pf Since $w$ has Abelian complexity $n$, it contains all $n$ letters. It follows that  $G$ must be connected. This implies that $G$ contains a spanning tree. The spanning tree contains $n-1$ edges. Since the factors of $w$ of length 2 represent exactly $n$ Abelian words, $G$ contains exactly $n$ edges. $\Box$

Let $b\in\Sigma$. Define $T(b)=\{[abc]:abc\mbox{ is a factor of }w\mbox{ for some }a,c\in\Sigma.\}$. We call an element of $T(b)$ a {\bf triple associated with} $b$. Since $w$ is recurrent, every letter of $\Sigma$ occurs in $w$ as the middle letter of at least one factor of length 3. This means that each letter of $\Sigma$ has at least one triple associated with it.
\begin{lemma}\label{distinct triples}Suppose that $a\ne b$ but triple $[abc]$ is associated with both $a$ and $b$. Then exactly one of he following occurs:
\begin{enumerate}
\item{$abc$ is a triangle in $G$}
\item{$a=c$ and $G$ contains the loop $aa$.}
\item{$b=c$ and $G$ contains the loop $bb$.}
\end{enumerate}
\end{lemma}
\pf Since $abc$ is associated to $b$, at least one of $abc$ and $cba$ is a factor of $w$. It follows that $ab$ and $bc$ are edges of $G$. Since $abc$ is a triple associated to $a$, at least one of $bac$ and $cab$ is
a factor of $w$, so that $ca$ is also an edge of $G$. If $a$, $b$ and $c$ are distinct, then $G$ contains  triangle $abc$. If two of them are the same, then since $a\ne b$, one of $bc$ and $ca$ is a loop.$\Box$

\begin{lemma}\label{degree 3} Suppose that $b$, $c$ and $d$ are distinct neighbours of $a$ in $G$. Then $|T(a)|\ge 2$.
\end{lemma}
\pf Since $b$ is a neighbour of $a$, then either $ba$ or $ab$ is a factor of $w$. By recurrence, either $bax$ or $xab$ will therefore be a factor of $w$ for some $x\in\Sigma$, and $[xab]$ is a triple associated with $a$. Similarly,  $a$ must have associated triples $[cay]$, $[daz]$ for some letters $y,z\in \Sigma$. If $[bax]\ne[cay]$, then we are done. Otherwise, $x=c$ and $y=b$ so that $[daz]\ne[bax].\Box$

\begin{lemma} \label{degree 2}Suppose that $a$ has distinct neighbours $b$ and $c$ in  $G$. Either $[bac]\in T(a)$ or $|T(a)|\ge 2$.
\end{lemma}
\pf One of $ab$ and $ba$ is a factor of $G$. Suppose $ab$ is a factor of $w$. (The other case is similar).  By recurrence, $w$ has a factor $xab$ for some $x$. If $x=c$, then $[bac]$ is associated to $a$, and we are done. Otherwise, $[xab]\in T(a)$, and $T(a)$ also includes a triple involving $c.\Box$

\subsection*{Case 1: Cycle $C$ is an $m$-cycle, $m\ge 4$.}
By Lemma~\ref{one cycle},$C$ is the unique cycle in $G$. This implies that $G$ contains no loops or triangles, so that triples associated with distinct vertices are distinct by Lemma~\ref{distinct triples}. At least one triple is associated with each of the $n$ vertices of $G$. Since the
Abelian complexity of $w$ is exactly $n$, the total number of triples associated with the vertices of $G$ is $n$. We conclude that $|T(a)|=1$ for each $a\in\Sigma$.  From Lemma~\ref{degree 3} we conclude that each vertex  of $C$ has degree exactly 2, so that $C$ is a connected component of $G$. Since $G$ is connected, $G=C$ is an $n$-cycle. Without loss of generality let the vertices be connected in the natural order $123\cdots n1$. By Lemma~\ref{degree 2}, we conclude that the triples associated with the vertices of $G$ are
$[123],[234],[345],\ldots,[(n-2)(n-1)n],[(n-1)n1],[n12]$. Since $w$ must be walked on $G$ respecting the possible triples, we conclude that $w$ is a suffix of $(123\cdots n)^\omega$ or of   $(n\cdots 321)^\omega$ and thus has period $n$. However, this means that $w$ contains exactly one factor of length $n$, up to anagrams. This is a contradiction.
\subsection*{Case 2: Cycle $C$ is a loop.}
Without loss of generality, let the loop edge be $11$. At least one triple is associated with each of the $n$ vertices of $G$.
A triple of the form $[111]$ could only ever be associated to $1$. Also, $b,c$ are neighbours of $1$ and $[b1c]$ is associated to $b$, then $1bc$ or $cb1$ is a factor of $w$. This implies that $bc$ is an edge of $G$ so that $ibc$ is a triangle (if $b\ne c$) or $bc$ is a loop (if $b=c$). Since $11$ is the only cycle in $G$ by Lemma~\ref{one cycle}, this is impossible. It follows that triples of the form $111$ or $b1c$ with $b$ and $c$ neighbours of $1$ can only ever be associated to $1$. It now follows that $1$ can be associated to at most a single triple of the form $111$ or $b1c$ where $b,c$ are neighbours of $1$; if $1$ is associated to two such triples $T_1$ and $T_2$, then each of the $n-1$ other vertices of $G$ is associated to a triple, and these triples are distinct from $T_1$ and $T_2$, and from each other by Lemma~\ref{distinct triples}. Then, however, we have at least $n+1$ distinct triples, violating the Abelian complexity of $w$.

We make cases based on whether $1$ is associated to a triple of the form $111$ or $b1c$ where $b$ and $c$ are neighbours of $1$.
\subsubsection*{Case 2a: Vertex $1$ is associated to a triple of the form $[111]$.}
Each vertex of $G-\{1\}$ is associated to some triple other than $[111]$, and those triples are distinct from each other and from $111$. Let $b$ be a neighbour of $1$. At least one of $b1$ and $1b$ is a factor of $w$. Since $1$ is not associated to any triple $[b1c]$ where $b$ and $c$ are neighbours of 1, it follows that $b11$ or $11b$ must be a factor of $w$. Since the Abelian complexity of $w$ is $n$, we conclude that $[11b]=[1b1]$ is the unique triple associated with $b$. From Lemma~\ref{degree 2}, it follows that $1$ is the only neighbour of $b$. Graph $G$ is therefore the star with center $1$; the edges of $G$ are precisely $E(G)=\{1k:1\le k\le n\}$.

Let $m$ be least such that $w$ has a factor $d1^me$ where $d,e\ne 1$. Without loss of generality, say that $21^m3$ is a factor of $w$. Let $b$ be any vertex of $G-\{1\}.$ Since $1b$ is an edge of $G$, $w$ has a factor $b1$ or $1b$, hence a factor $1^2b1^m$ or $1^mb1^2$. (Recall that 1 is the only neighbour of $b$ in $G$.) It follows that up to anagrams, the $n$ factors of $w$ of length $m+3$ are $121^m3$, $1^m21^2$, $1^m31^2$, $1^m41^2,\ldots$, $1^mn1^2$. In particular, $w$ has no factor $1^{m+2}$, and in any factor of the form $b1^kc$ with $b,c\ne 1$ and $k\le m$ we must have $\{b,c\}=\{2,3\}$ and $k= m$.

Now consider the shortest factor of $w$ containing a letter from $\{2,3\}$ and a letter from $\{4, 5, \ldots, n\}$. By our last remark, it must have the form $b1^kc$ or $c1^kb$ where $b\in\{2,3\}$, $c\in\{4, 5, \ldots, n\}$ and $k\ge m+1$. Since $w$ has no factor $1^{m+2}$, $k=m+1$, and we have found an $(n+1)^{st}$ Abelian factor of $w$. This is a contradiction.
\subsubsection*{Case 2b: Vertex $1$ is associated to exactly one triple of the form $[b1c]$ where $b$ and $c$ are neighbours of $1$.}
In this case, $[111]$ is not associated with $1$; i.e., $111$ is not a factor of $w$. Each vertex of $G-\{1\}$ is again associated to some triple other than $[b1c]$, and these triples are distinct from each other. Let $d$ be a neighbour of $1$ other than $b$ or $c$. Vertex $1$ cannot be associated to another triple $[d1e]$ where $e$ is a neighbour of $1$. Therefore, at least one of $d11$ and $11d$ is a factor of $w$. It follows that $[11d]=[1d1]$ is the unique triple associated with $d$. We conclude that $1$ is the only neighbour of $d$; viz., $d$ is a leaf. We see also that (except possibly once at the beginning of $w$) $d$ always appears in $w$ in the context $11d11$.

Now consider the shortest factor of $w$ containing a letter from $\{b,c\}$ and a neighbour of $1$ other than $b$ or $c$. By our last remark, and relabeling $b$ and $c$ if necessary, this factor must have the form $c1^kd$ or $d1^kc$, $k\ge 2$, $d$ some neighbour of $1$ other than $b$ and $c$. If $[11c]$ is not associated to $c$, then $[11c]$ and $[b1c]$ are associated only to 1, and we can count $n+1$ distinct triples associated to vertices of $G$. This violates the Abelian complexity of $G$. We conclude that $[11c]$ must be the unique triple associated with $c$, and $c$ is a leaf.
\subsubsection*{Case 2bi: Vertex $b$ is a leaf.}
In this case, each neighbour of 1 is a leaf; graph $G$ is the star with center $1$. The edges of $G$ are precisely $E(G)=\{1k:1\le k\le n\}$.
Let $m$ be least such that $w$ has a factor $x1^md$ or $d1^mx$ where $x \in\{b,c\}, d\notin\{1,b,c\}$. Since $[b1c]$ is the unique triple associated only with 1, $m\ge 2$. On the other hand $111$ is not a factor of $w$, so $m=2$. Without loss of generality, assume that $b=2$, $c\ne 3$, and $2113$ or $3112$ is a factor of $w$. Since $1$ is the only neighbour of $2$, it follows that $12113$ or $31121$ is a factor of $w$. We have already seen that $11d11$ is a factor of $w$ if $d\ne 1, b, c$. It follows that, up to anagrams, the following $n+1$ factors of length 4 appear in $w$:
$$121c, 2113, 1121, 1131, 1141, \ldots, 11n1.$$
This is a contradiction.
\subsubsection*{Case 2bii: Vertex $b$ has degree at least 2.}
Since our Abelian complexity is $n$, and triple $b1c$ is associated only with $1$, for any vertex $d$ of $G-\{1\}$, $|T(d)|=1$. By Lemma~\ref{degree 3}, deg$(d)\le 2$. We may therefore assume that the edges of $G$ are
$$11, 1n, 1(n-1), 1(n-2), \ldots 1(r+1), 1r, r(r-1), (r-1)(r-2), (r-2)(r-3),\ldots, 32$$
and the triples are $$[11n],[11(n-1)],[11(n-2)],\ldots$$ $$[11(r+1)],(r+1)1r],[1r(r-1)],[r(r-1)(r-2)],\ldots,[432],[323].$$
(We have $c=r+1$, $b=r$.)
It follows that up to anagrams, $w$ has the $n$ length 4 factors
$$\begin{array}{ll}11n1, 11(n-1)1,\ldots, 11(r+1)1&\mbox{ ($n-r$ factors)}\\
1(r+1)1r, (r+1)1r(r-1)&\mbox{ (2 factors)}\\
1r(r-1)(r-2), r(r-1)(r-2)(r-3),\ldots, 5432, 4323&\mbox{ ($r-2$ factors)}.
\end{array}$$
Now however, consider the shortest factor of $w$ containing letters from both $\{n, n-1,\ldots, r-2\}$ and $\{r+1,r\}$. This must have the form $x1^ky$ or $y1^kx$ where $x\in \{n, n-1,\ldots, r-2\}$ and $y\in\{r+1,r\}$. Since $[b1c]$ is the unique triple associated only to $1$, we cannot have $k=1$. Since $111$ is not a factor of $w$, we must have $k=2$. This gives an $(n+1)^{st}$ length 4 Abelian word in $w$, namely $x11y$. This is a contradiction.
\subsubsection*{Case 2c: Every triple associated with 1 has the form $[11b]$, $b\ne 1$.}
In this case, $w$ has no factors 111 or $b1c$ with $b,c\ne 1$. If $b$ is any neighbour of $1$ therefore, either $b11$ or $11b$ is a factor of $w$.
If $1$ has no non-leaf neighbour, $G$ is a star centered at $1$; for $2\le k\le n$, the only length three factors of $w$ containing $k$ are among $11k$, $1k1$ and $11k$. The triples of $G$ are precisely those of the form $[11k]$, $k\ne 1$, and $G$ has only $n-1$ triples. This is a contradiction.

Therefore, let $b$ be a non-leaf neighbour of 1, let $b'\ne 1$ be a neighbour of $b$. The shortest factor of $w$ containing $1$ and $b'$ must be $1bb'$ or $b'b1$, so $[1bb']$ is a triple associated with $b$. Now every vertex of $G-\{1\}$ has at least one triple associated with it, and all such associated triples must be distinct. Moreover, $b$ has triples $[11b]$ and $[1bb']$ associated with it. We have now listed $n$ distinct triples associated with the vertices of $G-\{1\}$. If 1 had another non-leaf neighbour $c\ne b$, then an $(n+1)^{st}$ triple $[1cc']$ would be associated to $c$. Since this is impossible,
it follows that $b$ is the only non-leaf neighbour of $1$.

Without loss of generality, let the neighbours of $1$ be exactly $2,3,\ldots r=b$, and let $r+1$ be a neighbour of $r$. The $r$ triples $[112],[113],\ldots[11r],[1r(r+1)]$ will be associated to vertices $1,2,\ldots, r$, while the triples associated with vertices $(r+1), (r+2),\ldots, n$ must be distinct from these and from each other. This means that exactly one triple is associated to each of vertices $(r+1), (r+2),\ldots, n$, so that by Lemma~\ref{degree 3}, they each have degree at most 2. Without loss of generality we may thus assume that
the edges of $G$ are
$$11, 12, 13, \ldots$$ $$1r, r(r+1), (r+1)(r+2), (r+2)(r+3),\ldots, (n-1)n$$
some $r$, $1<r\le n$. The $n$ triples associated to vertices of $G$ must thus be precisely
$$[112],[113],[114],$$ $$[11r],[1r(r+1)],[r(r+1)(r+2)],\ldots,[(n-2)(n-1)n], [(n-1)n(n-1)].$$
For $2\le k\le r-1$, the only neighbour of vertex $k$ is vertex $1$. Since $w$ has no factors 111 or $b1c$ with $b,c\ne 1$, it follows that
$w$ has a factor $11k1$ for $2\le k\le r-1$. In addition to these $r-2$ factors of length 4, the specification of triples forces $w$ to have (up to reversal) the $n-r+1$ factors $$11r(r+1),1r(r+1)(r+2),r(r+1)(r+2)(r+3),(r+1)(r+2)(r+3)(r+4),$$ $$\ldots, (n-3)(n-2)(n-1)n, (n-2)(n-1)n(n-1).$$
In addition, since $11r$ or $r11$ is a factor of $w$, so is a word $c11r$ or $r11c$, where $c$ is some neighbour of $1$ in $G$. This brings the count of Abelian factors of length 4 to $(r-2)+(n-r+1)+1=n$. Suppose now that $d$ is a neighbour of $1$ other than $r$ and $c$. Then $w$ contains a factor $d11$ or $11d$, hence a word $d11e$ or $e11d$, where $e$ is a neighbour of $1$. This brings the number of length 4 Abelian factors of $w$ to $n+1$, which is a contradiction. It follows that the only neighbours of 1 are $r$ and $c$. (Note that perhaps $r=c$.) The length 4 Abelian factors of $w$ are thus
$$ [c11r], [11r(r+1)], [1r(r+1)(r+2)],[r(r+1)(r+2)(r+3)],[(r+1)(r+2)(r+3)(r+4)],\ldots, $$ $$[(n-3)(n-2)(n-1)n], [(n-2)(n-1)n(n-1)].$$
In the case that $c\ne r$, this forces $w$ to be a suffix of $$\left(c11r(r+1)(r+2)(r+3)(r+4)\cdots\right.$$ $$ \left.(n-2)(n-1)n(n-1)(n-2)\cdots(r+4)(r+3)(r+2)(r+1)r11\right)^\omega,$$ and $w$ is periodic, with period $2n$. However, this means that $w$ contains exactly one factor of length $2n$, up to anagrams, which is a contradiction.\\
In the case that $c= r$, this forces $w$ to be a suffix of $$\left(r(r+1)(r+2)(r+3)(r+4)\cdots\right.$$ $$\left.(n-2)(n-1)n(n-1)(n-2)\cdots(r+4)(r+3)(r+2)(r+1)r11\right)^\omega,$$ and again $w$ is periodic, with a contradiction.
\subsection*{Case 3: Cycle $C$ is a $3$-cycle.}
Let the vertices of $C$ be $a, b, c$. By Lemma~\ref{distinct
triples}, the only triple which can be associated with more than
one vertex is $[abc]$. Each vertex of $G-\{a,b,c\}$ is associated
with some triple, and these must all be distinct. This accounts
for $n-3$ triples. Since $G$ is connected, $w$ must contain some
factor of the form $xab$, $xba$, $xbc$, $xcb$, $xca$ or $xac$, for
some $x\notin \{a,b,c\}$. Suppose without loss of generality that
$[xab]$ is associated to $a$ for some $x\not\in\{a,b,c\}$. Then $a$
has degree at least $3$, so that $|T(a)|\ge 2$ by
Lemma~\ref{degree 3}. Since $|T(b)|,|T(c)|\ge 1$ but the total
number of distinct triples associated to vertices of $G$ is $n$,
at least two of $a$, $b$ and $c$ have an associated triple in
common. That triple must be $[abc]$. So far, we have found that
$[xab],[abc]\in T(a)\cup T(b)\cup T(c).$

\subsubsection*{Case 3a: $T(a)\cup T(b)\cup
T(c)=\{[abc],[bax]\}$.} The shortest factor of $w$ starting with
$x$ and ending in one of $b$ or $c$ will be $xab$. (Such a factor exists
because $w$ is recurrent.) Let $uxab$ be a prefix of $w$. As the
only triple associated to $b$ is $[abc]$, $w$ has $uxabc$ as a
prefix. Again, $T(c)=\{[abc]\}$, so $uxabca$ is a prefix of $w$.
The only triple in $T(a)$ having $c$ as one of its letters is
$[abc]$, so $uxabcab$ is a prefix of $w$. Continuing in this way,
we find that $w=ux(abc)^\omega$. This is impossible, since $x$
must appear in $w$ infinitely often, by recurrence.

\subsubsection*{Case 3b: $|T(a)\cup T(b)\cup
T(c)|=3$.} The argument of {\bf Case 3a} can still be applied if
we add to $T(a)\cup T(b)\cup T(c)$ another triple from $bab$,
$cbc$ or $aca$; none of these triples allows us to break the
circular order $a-b-c-a$ on $\{a,b,c\}$ which commences with
$xab$. Similarly, adding to $T(a)\cup T(b)\cup T(c)$ a triple
$[ybc]$ where $y\ne\{a,c\}$ is a neighbour of $b$ would lead to
the same contradiction. Again a triple $[yca]$ where $y\ne\{b,a\}$
is a neighbour of $c$, or a triple $[yab]$ where $y\ne\{b,c\}$ is
a neighbour of $a$ leads to a contradiction.

We may therefore assume $w$ contains a factor $aba$, $bcb$, $cac$,
or a factor of the form $aby$, $bcy$ or $cay$, $y\not\in
\{a,b,c\}$.

Since $T(a)\cup T(b)\cup T(c)$ contains three distinct triples,
and a distinct triple is associated to each of the $n-3$ vertices
of $G-\{a,b,c\}$, we deduce that each vertex of $G-\{a,b,c\}$ is
only associated with a single triple, and thus has degree at most
2 by Lemma~\ref{degree 3}. We recall that $G$ contains exactly one
cycle. Graph $G$ therefore consists of the triangle $abc$,
together with 1 or more paths radiating from its vertices.

\subsubsection*{Case 3bi: Word $w$ contains a factor $aba$, $bcb$, $cac$.}

In this case, the only vertex of $G-\{a,b,c\}$ which is adjacent
to any of $a$, $b$ and $c$ is $x$. Graph $G$ consists of triangle
$abc$ together with a single path adjacent to $a$. Relabel
$x=a_1$, and let the edges of $G$ be
$$ab, bc, ca, aa_1, a_1a_2, a_2a_3,\ldots, a_{r-1}a_r$$
where $r=n-3$.\footnote{If $r=1$, we use the convention $a_0=a_2=a$.} By Lemma~\ref{degree 2}, the triples of $G$ are
precisely
$$[baa_1],[aa_1a_2],[a_1a_2a_3],[a_2a_3a_4],\ldots,[a_{r-2}a_{r-1}a_r],[a_{r-1}a_ra_{r-1}],[abc],T$$
where $T$ is one of $[aba]$, $[bcb]$ and $[cac].$

The only two triples containing both $a$ and $a_1$ are $[baa_1]$
and $[aa_1a_2]$. We conclude that $w$ contains Abelian factor
$[aba_1a_2]$. Reasoning similarly, we find that for $r\ge 4$, the following $n-2$
length 4 Abelian factors must be in $w:$
$$[baa_1a_2],[aa_1a_2a_3],[a_1a_2a_3a_4],\ldots,[a_{r-3}a_{r-2}a_{r-1}a_r],[a_{r-2}a_{r-1}a_ra_{r-1}].$$
(The stipulation $r\ge 4$ is only for notational convenience. If $r=3$, let $a_0=a$, and the 3 length 4 Abelian factors are $[baa_1a_2],[aa_1a_2a_3],[a_1a_2a_3a_2]$. If $r=2$, let $a_0=a$, $a_{-1}=b$, and the 2 length 4 Abelian factors are $[baa_1a_2],[aa_1a_2a_1].$ Finally, if $r=1$, a length 4 Abelian factor is $[baa_1a].$)

Now consider a factor $v$ of $w$ of the form $a_1\{a,b,c\}^*a_1$, containing $ac$ or $ca$ as a factor. Such a $v$ exists by recurrence.
 Since the only
triple joining $abc$ and $G-\{a,b,c\}$ is $[baa_1]$, $v$ can be
written in the form $a_1abv_1baa_1$ where  $v_1\in\{a,b,c\}^*$ and
$ac$ or $ca$ appears in $v_1$. The circular order of $abv_1ba$
changes exactly once, from $a-b-c-a$ to $a-c-b-a$, at triple $T$.
Thus $v_1$ cannot both begin and end with $a$, lest $aba$ appear
twice in $abv_2ba$. Thus $v_1$ must either begin or end with $c$,
so that $a_1abc$ or $cbaa_1$ is a factor of $w$, yielding Abelian
factor $[a_1abc]$ in either case. Notice that we have shown that
$abv_1ba$ cannot both begin and end with a
palindrome.\footnote{Throughout, when we say ``palindrome'' we mean
  one of the three palindromes $aba, bcb, cac$.} It thus follows that
$[tz]$ is also an Abelian factor of $w$, where $t\in T$ is a palindrome and $z$ is the letter of $\{a,b,c\}$ not appearing in $t$.

Suppose now that $v_1$ begins with $a$. The case where $v_1$ ends in $a$ is similar.
Then $a_1aba$ is a factor of $w$, and we have enumerated all $n$ length 4 Abelian factors of $w$: the $n-3$ previously listed, plus $[a_1abc],[tz]=[abac]$ and $[a_1aba]$. The circular order of $\{a,b,c\}$ in $abv_1ba$ changes exactly once (with $aba$), so that $abv_1ba\in aba(cba)^+$. It follows that $bacba$ is a suffix of $abv_1ba$, and $w$ also contains Abelian factor $[bacb]$. This is a contradiction. We conclude that $v_1$ cannot begin or end with $a$, and hence must begin and end with $c$. Since $[ca]$ is an Abelian factor of $v$, we cannot have $abv_1ba=abcba$.

Thus far, $w$ has Abelian factors $[a_1abc]$ and $[tz]$ in addition to the $n-3$ length 4 Abelian factors previously listed.
Let $y$ be the central letter of palindrome $t$ and write
 $abv_1ba=v_2yv_3$ where $v_2y$ is a prefix of $(abc)^\omega$ and $yv_3$ is a suffix of $(cba)^\omega$. We must have $|v_2|\equiv|v_3|$ (mod 3). Also, $|v_2|,|v_3|\ge 2$.  Suppose that $|v_2|>|v_3|$. Then $|v_2|\ge|v_3|+3\ge 5.$ In this case, $abv_1$ has a prefix $abcabc$, and $w$ contains Abelian factors $[abca]$, $[bcab]$, $[cabc]$. One of these is $[tz]$, but this still gives $n+1$ length 4 Abelian factors of $w$, which is impossible. We similarly rule out $|v_2|<|v_3|$. Note that we may also assume that $|v_2|\le 4$. Since $|abv_1ba|>5$, we find that $3\le |v_2|=|v_3|\le 4$. If $|v_2|=3$, then $abv_1ba=abcacba$, $t=cac$ and $[abca]$ and $[bcac]$ are Abelian factors of $w$. We have now specified all length 4 Abelian factors of $w$; none of these is $[a_1t]$, the central letter in $t$ is not $c$ and the set of length 4 Abelian factors of $w$ turns out to be determined by $t$ and $|v|=2|v_2|+3$. Similarly, if $|v_2|=4$, then $abv_1ba=abcabacba$, $t=aba$ and $[abca]$ and $[bcab]$ are Abelian factors of $w$. Again the set of all length 4 Abelian factors of $w$ is determined by $t$ and $|v|$, none of the Abelian factors is $[a_1t]$ and $c$ is not the central letter in $t$.

Since the two different possible lengths for $v$ give different sets of Abelian factors in $w$, it follows that $w$ contains exactly one factor $v$ of the form $a_1\{a,b,c\}^*a_1$ containing $[ac]$ as an Abelian factor. Now let $v'$ be any factor of $w$ of the form $a_1\{a,b,c\}^+a_1$. Word $v'$ must have prefix $a_1ab$ and suffix $baa_1$. However, since $[a_1t]$ is not an Abelian factor of $w$, $v'$ cannot have $a_1aba$ as a prefix or $abaa_1$ as a suffix. Word $v'$ therefore has $a_1abc$ as a prefix and $cbaa_1$ as a suffix. Again, $v'\ne a_1abcbaa_1$, since the central letter of $t$ is not $c$. We deduce that $v'$ has prefix $a_1abca$ or suffix $acbaa_1$, and must contain $[ac]$ as an Abelian factor. In summary, $w$ contains exactly one factor of the form $v=a_1\{a,b,c\}^*a_1$. If $r=1$, this shows that $w$ is periodic, giving a contradiction. If $r\ge 2$,
our earlier specification of the $n-3$ triples of $w$ along the path $baa_1\cdots a_r$ shows that $w$ contains a single factor of the form $a_1(\Sigma-\{a,b,c\})^*a_1$, namely $a_1a_2\cdots a_{r-1}ra_{r-1}\cdots a_2a_1$. Since $aa_1a$ is not a factor of $w$, we again deduce that $w$ is periodic, giving a contradiction.

\subsubsection*{Case 3bii: Word $w$ contains a factor $aby$, $bcy$ or $cay$, $y\not\in
\{a,b,c\}$.} We consider first the case where $w$ contains a
factor $aby$, $y\not\in \{a,b,c\}$. Since $b$, $c$ and $x$ are
neighbours of $a$, and $abc$ is the only cycle in $G$, we cannot
have $y=x$. Graph $G$ consists of triangle $abc$ together with
two paths adjacent to $a$ and $b$. Relabel $x=a_1$, $y=b_1$ and
let the edges of $G$ be
$$ab, bc, ca, aa_1, a_1a_2, a_2a_3,\ldots, a_{r-1}a_r, bb_1, b_1b_2, \ldots, b_{s-1}b_s$$
where $r+s=n-3$. By Lemma~\ref{degree 2}, the $n$ triples of $G$ are

$$[baa_1],[aa_1a_2],[a_1a_2a_3],[a_2a_3a_4],\ldots,[a_{r-2}a_{r-1}a_r],[a_{r-1}a_ra_{r-1}],[abc]$$
and
$$[abb_1],[bb_1b_2],[b_1b_2b_3],[b_2b_3b_4],\ldots,[b_{s-2}b_{b-1}b_s],[b_{s-1}b_sb_{s-1}].$$

The following $n-3$ length 4 Abelian factors must be in $w$:
$$[baa_1a_2],[aa_1a_2a_3],[a_1a_2a_3a_4],\ldots,[a_{r-3}a_{r-2}a_{r-1}a_r],[a_{r-2}a_{r-1}a_ra_{r-1}]$$
and
$$[abb_1b_2],[bb_1b_2b_3],[b_1b_2b_3b_4],\ldots,[b_{s-3}b_{s-2}b_{s-1}b_s],[b_{s-2}b_{s-1}b_sb_{s-1}].$$
Let $v$ be a shortest factor of $w$ of the form $\{a_1,b_1\}uc$. A
prefix of $v$ must be $a_1ab$ or $b_1ba$. Suppose $v$ has prefix
$a_1ab$. Then $v\in a_1ab\{a,b,c\}^*c$. Letters $a,b,c$ must have
circular order $a-b-c-a$ in $v$, since there are no palindromes in
$v$ to change the circular order, and $v$ starts $a_1abc$. Let $p$
be a prefix of $w$ of the form $qa_1(abc)^j$ with $j$ as large as
possible. If $j=2$, then $abcabc$ is a factor of $w$, so that $w$
contains 4 more Abelian factors:
$$[a_1abc],[abca],[bcab],[cabc].$$
This is impossible, since then $w$ has $n+1$ distinct length 4
Abelian factors. It follows that $j=1$, and $qa_1abcabb_1$ is a
prefix of $w$. (Recall that the only triples associated with
$a, b$ or $c$ are $[abc],[a_1ab],[abb_1]$.) Now, however, $w$
contains Abelian factors $$[a_1abc],[abca],[bcab],[cabb_1],$$
again giving a contradiction.

Now consider the case where $w$ contains a factor $bcy$, $y\not\in
\{a,b,c\}$. Since $abc$ is the only cycle in $G$, $y$ is not a
neighbour of $b$ or $a$. Since $ac$ is an edge of $G$, either $ac$
or $ca$ is a factor of $w$. Suppose $ac$ is a factor of $w$. (The
other case is similar.)  Recall that the only triples associated
with one of $a$, $b$ or $c$ are $[abc],[xab]$ and $[bcy]$. The
only one of these containing both $a$ and $c$ is $[abc]$. If $ac$
is a factor of $w$, then it must therefore be preceded and
followed by $b$, and occurs in the context $bacb$. Since neither
of $[xab]$ and $[bcy]$ is associated with $b$, $cb$ is followed by
$a$. Again. $ba$ is preceded by $c$, so $ac$ occurs in the context
$cbacba$. Since $w$ is recurrent, it cannot have $(cba)^\omega$ as
a suffix. It follows that $w$ must have a factor $cbacbax$. This
implies that $[cbac],[bacb],[acba],[cbax]$ are length 4 Abelian
factors of $w$. As in previous cases, the paths attached to
vertices $a$ and $c$ of triangle $abc$ furnish another $n-3$
distinct length 4 Abelian factors, giving a contradiction.

The final case occurs when $w$ contains a factor $cay$, $y\not\in
\{a,b,c\}$. In this case $G$ consists of a triangle with two
disjoint paths attached at $a$. In the usual way, we find $n-3$
length 4 Abelian factors of $w$, each containing at least two
path vertices (i.e. vertices of $G-\{a,b,c\}$). If $abcabc$ or
$cbacba$ were a factor of $w$, $w$ would then contain 4 additional
length 4 Abelian factors $[abca]$, $[bcab]$, $[cabc]$ and
$[bcay]$, giving a contradiction. We therefore conclude that the
only factor of $w$ of the form $x\{a,b,c\}^*y$ is $xabcay$, and the
only factor of $w$ of the form $y\{a,b,c\}^*x$ is $xabcay$. Thus,
if $L_1$ is the leaf at the end of the path starting with $a-x$
and $L_2$ is the leaf at the end of the path starting with $a-y$,
$w$ has only one factor of the form
$L_1(\Sigma-\{L_1,L_2\}^*)L_2$, and only one factor of the form
$L_2(\Sigma-\{L_1,L_2\}^*)L_1$, so that $w$ is periodic,
oscillating between $L_1$ and $L_2$. The
periodicity of $w$ gives a contradiction.$\Box$

\end{document}